\documentclass[11pt,twoside]{article}
\usepackage{a4}
\usepackage{amssymb,amsmath,amsthm,latexsym}
\usepackage{amsfonts}
\usepackage{amsfonts}
\usepackage{graphicx}
\usepackage[utf8]{inputenc}  
\usepackage[T1]{fontenc}  

\usepackage{enumitem}

\newcommand{\norm}[1]{\lVert #1 \rVert}
\usepackage{tikz-cd}

\newtheorem{theorem}{Theorem}[section]

\newtheorem{corollary}[theorem] {Corollary}
\newtheorem{definition}[theorem]{Definition}
\newtheorem{example}[theorem]{Example}
\newtheorem{lemma} [theorem]{Lemma}

\newtheorem{problem}[theorem]{Problem}
\newtheorem{proposition}[theorem]{Proposition}
\newtheorem{remark}[theorem]{Remark}

\vspace{5cm}

\begin{document}
	
	\label{'ubf'}  
	\setcounter{page}{1}                                 

	\markboth {\hspace*{-9mm} \centerline{\footnotesize \sc
			  }
	}
	{ \centerline                           {\footnotesize \sc  
Qi Liu,	Yuxin Wang, Mengmeng Bao                                                } \hspace*{-9mm}              
	}

	\vspace*{-2cm}

	\begin{center}
		{ 
			{\Large \textbf { \sc  How orthogonality influences geometric constants?
				}
			}
			\\

			\medskip

			{\sc Qi Liu }\\
			{\footnotesize School of Mathematics and Statistics, Anqing Normal University, Anqing 246133, China}\\
			{\footnotesize e-mail: {\it liuq67@aqnu.edu.cn}}
			
			{\sc  Yuxin Wang }\\
			{\footnotesize School of Mathematics and Statistics, Anqing Normal University, Anqing 246133, China}\\
			{\footnotesize e-mail: {\it y24060028@stu.aqnu.edu.cn}}
			
			{\sc Mengmeng Bao }\\
			{\footnotesize School of Mathematics and Statistics, Anqing Normal University, Anqing 246133, China}\\
			{\footnotesize e-mail: {\it y25060053@stu.aqnu.edu.cn}}
		}

	\end{center}

	\thispagestyle{empty}

	\hrulefill

	\begin{abstract}  
		{\footnotesize In this paper, based on isosceles orthogonality, we have found equivalent definitions for four constants:  $A_2(X)$  proposed by Baronti in 2000 [J. Math. Anal. Appl. 252(2000), 124-146],  $C'_{\mathrm{NJ}}(X)$  introduced by Alonso et al. in 2008 [Stud. Math. 188(2008), 135-150], $T(X)$ introduced by Alonso et al. in 2008 [J. Math. Anal. Appl. 340(2008), 1271-1283] and  $L'_{\mathrm{YJ}}(X)$  put forward by Liu et al. in 2022 [Bull. Malays. Math. Sci. Soc., 45(2022), 307-321]. A core commonality among these four constants is that they
			are all restricted to the unit sphere. This finding provides us with the
			following insight: could it be that several constants defined over the whole
			space, when combined with suitable orthogonality conditions, are equivalent
			to their restrictions to the unit sphere?  Motivated by this question, we further study the corresponding problem for Birkhoff-James orthogonality. Because this orthogonality is generally non-symmetric, a direct replacement of isosceles orthogonality is not possible. We therefore introduce a norming-functional rectification method, which represents unit-sphere configurations by rectified Birkhoff-James orthogonal data. Consequently, exact Birkhoff-James orthogonal representations are obtained for the above four constants.
		}
	\end{abstract}
	\hrulefill



\section{Introduction}\label{sec1}

Geometric constants are crucial for functional analysis problems (see \cite{15,16,17,38,53}), and their connection to functional equation stability is emphasized in \cite{36}. In addition to classical geometric constants, the author \cite{37} introduced a geometric constant without infimum or supremum, which is very interesting. Apart from  this, scholars have gradually turned their attention to skew geometric constants. The parameter asymmetry makes skew geometric constants comparatively more difficult to study, and relevant research findings can currently be found in references \cite{02,46,01,47,45,43,44}.

In recent years, the academic community has defined and studied a large number of constants. Numerous proofs have been made regarding the properties and relationships of these constants, and the inequalities between them (sometimes extremely complex) have also been clearly explained. The core goal of this paper is to explore the role of isosceles orthogonality in Banach spaces. By proposing three isosceles orthogonal constants equivalent to the existing ones, we believe these new findings can help scholars re-examine the issue of orthogonality in Banach spaces.

We revisit some types of orthogonality concepts originally defined in normed linear spaces. In 1945, James \cite{03} first put forward the concept of isosceles orthogonality, denoted by $x \perp_I y$, which holds if and only if $\|x + y\| = \|x - y\|$. Inspired by the classical Pythagorean theorem, another orthogonal relation in a normed space $(X, \|\cdot\|)$, called Pythagorean orthogonality \cite{03}, can be defined as follows: $x \perp_P y$ when $\|x - y\|^2 = \|x\|^2 + \|y\|^2$. Although the definitions are different, these two types of orthogonality are equivalent in inner product spaces. Furthermore, in his work, Birkhoff \cite{05} defined a type of orthogonality known as Birkhoff-James orthogonality, which is defined as: for elements  $x, y$  in a normed linear space  $X$ ,  $x$  is Birkhoff-James orthogonal to  $y$  (denoted  $x \perp_B y$ ) if  $\|x + ty\| \geq \|x\|$  for all $t\in\mathbb{R}$, the authors \cite{40,42} introduce some properties and applications of Birkhoff-James orthogonality. In 1934, Robert \cite{41} proposed the concept of Robert orthogonality, defined as follows:  Two vectors \( x \) and \( y \) are Robert orthogonal if and only if  
$ \| x + \lambda y \| = \| x - \lambda y\| $  
for all \( \lambda\in \mathbb{R} \). In fact, Robert orthogonality is a generalization of isosceles orthogonality. For more definitions of orthogonality in normed linear spaces, readers may refer to works such as \cite{04,05,07,08,06} and the bibliographies therein. 

In Clarkson's work \cite{31}, he introduced the von Neumann-Jordan constant $C_{\mathrm{NJ}} (X )$:
$$C_{\mathrm{NJ}}(X)=\sup\left\{\frac{\|x+y\|^2+\|x-y\|^2}{2\|x\|^2+2\|y\|^2}:x,y\in X,(x,y)\neq(0,0)\right\}.$$
The definition of the modified von Neumann-Jordan constant is given below (Alonso et al. \cite{18}):
$$C_{\mathrm{NJ}}^{\prime}(X)=\sup\left\{\frac{\|x+y\|^2+\|x-y\|^2}{4}:x,y\in S_X\right\}.$$
Next, we give the definition of constant 
$A_2(X)$, which is closely related to the core proof of this paper. In their work \cite{21}, M. Baronti, E. Casini, and P. L. Papini introduced the constant  $A_2(X)$.
There are many interesting results concerning this constant, and for specific details, please refer to the relevant literature \cite{19}.
$$
A_2(X)=\sup \left\{\frac{\|x+y\|+\|x-y\|}{2}: x, y \in S_X\right\} .
$$
The authors \cite{53}  considers the triangles with vertices \(x\), \(-x\) and \(y\), where \(x, y\in S_X\).  By utilizing the geometric means of the variable lengths of the sides of these triangles, the following constant for Banach spaces is defined:
$$T(X)=\sup\limits_{x,y\in S_X}\bigg(\|x+y\|\|x-y\|\bigg)^\frac{1}{2}.$$
Moslehian and Rassias \cite{47} established a new equivalent characterization of inner product spaces, namely the Euler-Lagrange type identity, which replaces the parallelogram law: A normed space ($X,\|\cdot\|$) is an inner product space if and only if
$$
\|\tau x+\upsilon y\|^2+\|\upsilon x-\tau y\|^2=\left(\tau^2+\upsilon^2\right)\left(\|x\|^2+\|y\|^2\right)
$$
holds for all nonnegative real numbers $\tau, \upsilon$ and all $x, y \in X$. Inspired by this, Liu  et al. \cite{10}  generalized von Neumann-Jordan constant and introduced the following skew form of the von Neumann-Jordan constant:
$$
\begin{aligned}
	L_{\mathrm{YJ}}^{\prime}(\tau, \upsilon, X) & =\sup \left\{\frac{\|\tau x+v y\|^2+\|v x-\tau y\|^2}{2\left(\tau^2+v^2\right)}: x, y \in S_X\right\},
\end{aligned}
$$where $\lambda, \mu>0$.
We now summarize several properties of these constants. For a Banach space $X$, the following hold:

(i) The inequality $1 \leq L_{\mathrm{YJ}}^{\prime}(\tau, \upsilon, X) \leq 1+\frac{2 \tau \upsilon}{\tau^2+\upsilon^2}$ is satisfied;

(ii) $X$ is a Hilbert space if and only if $L_{\mathrm{YJ}}^{\prime}(\tau, \upsilon, X)=1$;

(iii) $X$ is uniformly non-square if and only if $L_{\mathrm{YJ}}^{\prime}(\tau, \upsilon, X)<1+\frac{2 \tau  \upsilon}{\tau^2+ \upsilon^2}$.

It is noteworthy that Yang et al. determined the exact value of this constant in the regular octagon space in reference \cite{09}.

In \cite{23}, Joly introduced a constant based on Birkhoff-James orthogonality named rectangular constant. In a more recent work \cite{22}, Baronti, Casini, and Papini have put forward a new constant referred to as the isosceles constant. This definition bears a strong resemblance to that of the rectangular constant, except that it utilizes isosceles orthogonality rather than Birkhoff orthogonality.
$$
H(X)=\sup \left\{\frac{1+\lambda}{\|x+\lambda y\|}: x, y \in S_X, x \perp_I y, \lambda \geq 0\right\}.
$$

Orthogonality plays a crucial role in the research of Banach space theory,
in addition to the isosceles orthogonality constants introduced above, scholars have proposed several other types of isosceles orthogonality constants in \cite{12,08,51,34}. This has provided us with numerous ideas for studying isosceles orthogonality constants:
$$D(X)=\inf\left\{\inf_{\lambda\in\mathbb{R}}\|x+\lambda y\|: x,y\in S_X,x\perp_{I}y\right\},$$
$$DW_I(X)=\sup\left\{\frac{\|x\|+\|y\|}{\|x-y\|}\left\|\frac{x}{\|x\|}-\frac{y}{\|y\|}\right\|:x,y\in X\backslash\{0\},x\perp_Iy\right\},$$
$$IB(X)=\inf\left\{\frac{\inf\limits_{\lambda\in\mathbb{R}}\|x+\lambda y\|}{\|x\|}:x,y\in X\backslash\{0\},x\perp_{I}y\right\},$$
$$L_X(\lambda)=\sup\left\{\frac{\|\lambda x+(1-\lambda)y\|^2+\|(1-\lambda)x+\lambda y\|^2}{\|x+y\|^2}:x,y\in X\backslash\{0\},x\perp_Iy\right\},$$where $\lambda\in[0,\frac12)$. But it is not limited to this. For example, isosceles orthogonality can also play a role in non-commutative geometry theory \cite{32}.

Through the research on the above constants, we have come up with the following idea: Can classical geometric constants that do not involve orthogonality conditions be expressed in terms of geometric constants that incorporate orthogonality conditions? In fact, existing research \cite{34,35} has indicated that both the von Neumann-Jordan constant and the $p$-th von Neumann-Jordan constant can be represented by isosceles orthogonal constants. Nevertheless, this paper establishes a connection between certain constants restricted to the unit sphere and isosceles orthogonal constants, thereby offering new perspectives on this issue. It will be concluded in this paper that the constants confined to the unit sphere are closely associated with the condition of isosceles orthogonality, which undoubtedly further corroborates the fact that a specific category of geometric constants is intimately linked to isosceles orthogonality.

\section{ The $\widetilde{H}(X)$ ,   $\widetilde{H}^2(X)$ and $\widetilde{T}(X)$  constants}
In this section, we introduce three new constants $\widetilde{H}(X)$, $\widetilde{H}^2(X)$ and $\widetilde{T}(X)$ based on isosceles orthogonality. Through research, we have found that the $\widetilde{H}(X)$ is equivalent to the classical constant $A_2(X)$, the  $\widetilde{H}^2(X)$ constant is equivalent to the  modified von Neumann-Jordan constant $C_{\mathrm{N J}}^{\prime}(X)$ and the $\widetilde{T}(X)$ is equivalent to the classical constant $T(X)$. In addition, building on previous studies, we have compared the relationships between these constants and some classical constants.

A slight modification of the constant $H(X)$, that is, adding the same parameter $\lambda$ to the orthogonal condition, leads to a new constant in the following form:
\begin{definition}
	$$\begin{aligned}
		\widetilde{H}(X)&=\sup \left\{\frac{\|x\|+\|\lambda y\|}{\|x+\lambda y\|}:x,y\in S_X, x \perp_I \lambda y,  \lambda \geq 0\right\}
		\\&=\sup \left\{\frac{1+\lambda }{\|x+\lambda y\|}: x,y\in S_X,x \perp_I \lambda y, \lambda \geq 0\right\}
		\\&=\sup \left\{\frac{\|x\|+\|y\|}{\|x+y\|}:x,y\in X,(x,y)\neq(0,0), x \perp_Iy \right\}.
	\end{aligned}	$$
\end{definition}

Clearly, the geometric constant $\widetilde{H}(X)$ also has the following equivalent definitions:
$$\widetilde{H}(X)=\sup \left\{\frac{2(\|x\|+\|y\|)}{\|x+y\|+\|x-y\|}:x, y \in X,(x,y)\neq(0,0), x \perp_I y\right\}.$$

\begin{remark}
	Since isosceles orthogonality does not possess the properties of Birkhoff-James orthogonality, it follows that $H(X)\neq\widetilde{H}(X)$. However,	in inner product spaces, since isosceles orthogonality satisfies  homogeneity, it can be concluded that the constant $H(X)$ and the constant $\widetilde{H}(X)$ are equivalent.
\end{remark}

\begin{proposition}If $X$ is a Banach space, then $\sqrt 2 \leq \widetilde{H}(X) \leq 2$.
\end{proposition}
\begin{proof}
	Since there exist $x, y \in S_X$ such that $\|x+y\|=\|x-y\|=\sqrt{2}$, we can find $x$ and $y$ satisfying $x \perp_I y$ and
	$
	\|x+y\|=\|x-y\|=\sqrt{2}.
	$
	We have 
	$
	\widetilde{H}(X)\geq\sqrt 2.
	$

	On the one hand, let $x, y \in X$ such that $x \perp_I y$. Put $$u=\frac{x+y}{\|x+y\|}, \quad v=\frac{x-y}{\|x+y\|},$$  then $\|u\|=\|v\|=1$, and we have $$u+v=\frac{2 x}{\|x+y\|}, u-v=\frac{2 y}{\|x+y\|}. $$
	Thus, 
	$$
	\begin{aligned}\frac{\|x\|+\|y\|}{\|x+y\|}&=\frac{1}{2}\left(\|u+v\|+\|u-v\|\right)
		\\& \leq 2,
	\end{aligned}
	$$
	which implies that $\widetilde{H}(X)  \leq 2$.
	
\end{proof}

The following result shows that $	\widetilde{H}(X)$ is equivalent to the classical constant  $A_2(X)$. This is an interesting result, indicating that isosceles orthogonality plays an important role in geometric constants and influences the equivalent forms of constants.

\begin{theorem}\label{tt}
	Let $X$ be a Banach space.  Then $\widetilde{H}(X)=A_2(X)$.
\end{theorem}
\begin{proof}
	First, for any $x,y\in X$ such that $x \perp_I y$, let
	$$
	u=\frac{x+y}{2}, v=\frac{x-y}{2}.$$
	We can deduce that
	$$
	\frac{\|x\|+\|y\|}{\|x+y\|}=\frac{\|u+v\|+\|u-v\|}{2\|u\|}
	$$
	and $\|u\|=\|v\|$. Let $$x^{\prime}=\frac{u}{\|u\|}, y^{\prime}=\frac{v}{\|v\|},$$ we have
	$$
	\begin{aligned}	\frac{\|x\|+\|y\|}{\|x+y\|}&=\frac{\|u+v\|+\|u-v\|}{2\|u\|}\\&=\frac{\left\|x^{\prime}+y^{\prime}\right\|+\left\|x^{\prime}-y^{\prime}\right\|}{2} \\&\leq  A_2(X).
	\end{aligned}$$
	which implies that
	$$
	\widetilde{H}(X) \leq A_2(X).
	$$
	On the other hand, for any $x, y \in S_X$, let  $u=\frac{x+y}{2}, v=\frac{x-y}{2}$, and hence $ \|u+v\|=\|u-v\|=1$. Furthermore, we can get that
	$$
	\begin{aligned}
		\frac{\|x+y\|+\|x-y\|}{2}&=\frac{2\|u\|+2\|v\|}{2 \|u+v\|}\\&=\frac{\|u\|+\|v\|}{\|u+v\|}\\& \leq \widetilde{H}(X).
	\end{aligned}
	$$
	Then $A_2(X) \leq \widetilde{H}(X) $, as desired.
\end{proof}

We recall that the Clarkson's modulus of convexity  \cite{26} for a space \( X \) is defined, for any \( \varepsilon \in [0, 2] \), as 
$$
\delta_X(\varepsilon) = \inf \left\{1 - \frac{\|x + y\|}{2} : x, y \in B_X,\ \|x - y\| \geq \varepsilon\right\},
$$
wherein "\( B_X \)" and "\( \geq \)" may equivalently be substituted with "\( S_X \)" and "\( = \)", respectively.

\begin{remark}
	In reference \cite{21}, it holds that $$A_2(X)=1+\sup \left\{\frac{\varepsilon}{2}-\delta_X(\varepsilon): \sqrt{2} \leq \varepsilon<2\right\},$$ and hence  we can deduce that 
	$$\widetilde{H}(X)=1+\sup \left\{\frac{\varepsilon}{2}-\delta_X(\varepsilon): \sqrt{2} \leq \varepsilon<2\right\}.$$  This indicates that the  constant  $\delta_X(\varepsilon)$ is also closely related to isosceles orthogonality.
\end{remark}

Subsequently, inspired by the $\widetilde{H}(X)$ constant mentioned above, we wondered whether we could define an isosceles orthogonal constant equivalent to the modified von Neumann-Jordan constant. Considering the relationship between  constant  $A_2(X)$  and constant  	$C_{\mathrm{N J}}^{\prime}(X) $, we defined the following isosceles orthogonal constant $	\widetilde{H}^{2}(X)$. Through research, it was found that this isosceles orthogonal constant is indeed equivalent to the modified von Neumann-Jordan constant.
\begin{definition}
	$$
	\widetilde{H}^{2}(X)=\sup \left\{\frac{\|x\|^2+\|y\|^2}{\|x+y\|^2}: x,y\in X,(x,y)\neq(0,0), x \perp_{I} y \right\}.
	$$
\end{definition}
It is not difficult to find that the above constant can measure the difference between isosceles orthogonality and Pythagorean orthogonality.

The following equivalent definition of the constant can be viewed as characterizing when the parallelogram law holds when elements in the space satisfy the isosceles orthogonality condition.

\begin{definition}\label{d2}
	$$
	\widetilde{H}^{2}(X)=\sup \left\{\frac{2\|x\|^2+2\|y\|^2}{\|x+y\|^2+\|x-y\|^2}:x,y\in X,(x,y)\neq(0,0), x \perp_{I} y\right\}.
	$$
\end{definition}

\begin{remark}
	It is not necessarily the case that 
	$$
	\inf \left\{\frac{2\|x\|^2+2\|y\|^2}{\|x+y\|^2+\|x-y\|^2}:x,y\in X, (x,y)\neq(0,0), x \perp_{I} y\right\}
	$$
	is always equal  $\frac{1}{\widetilde{H}^{2}(X)}$.
\end{remark}

Next, we will demonstrate that this definition is equivalent to the familiar von Neumann-Jordan constant restricted to the unit sphere. According to the method of Theorem \ref{tt},  we can obtain the following results.
\begin{theorem}\label{mt}
	Let $X$ be a Banach space.  Then $\widetilde{H}^2(X)=C_{\mathrm{N J}}^{\prime}(X)$.
\end{theorem}
\begin{proof}
	First, for any $x,y\in X$ such that $x \perp_I y$, let
	$$
	u=\frac{x+y}{2}, v=\frac{x-y}{2}.$$
	We can deduce that
	$$
	\frac{\|x\|^2+\|y\|^2}{\|x+y\|^2}=\frac{\|u+v\|^2+\|u-v\|^2}{4\|u\|^2}
	$$
	and $\|u\|=\|v\|$. Let $$x^{\prime}=\frac{u}{\|u\|}, y^{\prime}=\frac{v}{\|v\|},$$ we have
	$$
	\begin{aligned}	\frac{\|x\|^2+\|y\|^2}{\|x+y\|^2}&=\frac{\|u+v\|^2+\|u-v\|^2}{4\|u\|^2}\\&=\frac{\left\|x^{\prime}+y^{\prime}\right\|^2+\left\|x^{\prime}-y^{\prime}\right\|^2}{4} \\&\leq  C_{\mathrm{N J}}^{\prime}(X),
	\end{aligned}$$
	which implies that
	$$
	\widetilde{H}^2(X) \leq C_{\mathrm{N J}}^{\prime}(X).
	$$
	On the other hand, for any $x, y \in S_X$, let  $u=\frac{x+y}{2}, v=\frac{x-y}{2}$, and hence $ \|u+v\|=\|u-v\|=1$. Furthermore, we can get that
	$$
	\begin{aligned}
		\frac{\|x+y\|^2+\|x-y\|^2}{4}&=\frac{4\|u\|^2+4\|v\|^2}{4 \|u+v\|^2}\\&=\frac{\|u\|^2+\|v\|^2}{\|u+v\|^2}\\& \leq \widetilde{H}^2(X).
	\end{aligned}
	$$
	Then $C_{\mathrm{N J}}^{\prime}(X) \leq \widetilde{H}(X) $, as desired.
\end{proof}

\begin{theorem}
	$\widetilde{H}^2(X)=1$ if and only if $X$ is a Hilbert space.
\end{theorem}
\begin{proof}Since $\widetilde{H}^2(X)=1$, let $x, y \in S_X$, then $x+y \perp_I x-y$. We have
	$$
	\|x+y\|^2+\|x-y\|^2 \leq 4
	$$
	for any $x, y \in S_X$, and hence $X$ is a Hilbert space (see [1](page 48)).
	
	Conversely, it obviously holds.
\end{proof}

\begin{remark}
	The following approach to characterize inner product spaces is described in detail in the literature \cite{04}.
	$$
	x, y \in S_X, \quad x \perp_I \lambda y \Rightarrow\|x+\lambda y\|^2\sim1+\lambda^2
	$$
	characterizes inner product spaces, where "$\sim$" means either "$\leq$" or "$\geq$".
	We use the following equivalent form of the  constant $	\widetilde{H}^2(X)$: 
	$$
	\begin{aligned}
		\widetilde{H}^2(X) =\sup \left\{\frac{1+\lambda^2}{\|x+\lambda y\|^2}: x, y \in S_X,x \perp_I \lambda y,  \lambda \geq 0\right\}.
	\end{aligned}
	$$
	If $\widetilde{H}^2(X)=1$,  we have
	$
	\frac{1+\lambda^2}{\|x+\lambda y\|^2}\leq 1
	$ for any  $ \lambda \geq 0$
	where  $
	x \perp_I \lambda y,~ x, y \in S_X
	$.
	According to the above method, the constant $	\widetilde{H}^2(X)$ can similarly be used to characterize the inner product space.
	
\end{remark}

\begin{theorem}
	Let $X$ be a Banach space. Then, $\widetilde{H}^2(X)\leq C_{\mathrm{NJ}}(X)$.
\end{theorem}
\begin{proof}
	First, in \cite{25}, the authors have showed that $C_{\mathrm{NJ}}(X)$ can be written in the following equivalent form:
	
	$$
	C_{\mathrm{NJ}}(X)=\sup \left\{\frac{2\|x\|^2+2\|y\|^2}{\|x+y\|^2+\|x-y\|^2}:x,y\in X,(x, y) \neq(0,0)\right\}.
	$$
	For $x \perp_I y$, we have the following estimation:
	$$
	\begin{aligned}
		\frac{\|x\|^2+\|y\|^2}{\|x+y\|^2} 
		& = \frac{2\|x\|^2+2\|y\|^2}{\|x+y\|^2+\|x-y\|^2}
		\\&\leq C_{\mathrm{NJ}}(X).
	\end{aligned}
	$$
\end{proof}

\begin{remark}
	Building upon the classic von Neumann-Jordan constant, the authors \cite{24} introduced an orthogonal condition to propose this new constant concept. This innovative approach has provided significant inspiration for our subsequent introduction of the equivalent isosceles orthogonal constant and opened up new perspectives for related research directions.
	$$
	\mathrm{C}_{\mathrm{NJ}}^{I}(X)=\sup \left\{\frac{\|x+y\|^2+\|x-y\|^2}{2\left(\|x\|^2+\|y\|^2\right)}:x,y \in X, (x, y) \neq(0,0), x \perp_{I} y\right\}
	.$$Obviously, under the constraint of the isosceles orthogonal condition, the equivalence of $\mathrm{C}_{\mathrm{NJ}}^{I}(X)$ and $\widetilde{H}^2(X)$ does not hold.  However, it remains unknown whether the relation $\widetilde{H}^2(X)< \mathrm{C}_{\mathrm{NJ}}^{I}(X)$ also holds for this constant.
\end{remark}

It is noted that the constant $\widetilde{H}^2(X)$ can be expressed by the following equivalent definition:
$$
\widetilde{H}^2(X)=\sup \left\{\frac{2\|x\|^2+2\|y\|^2}{\|x+y\|^2+\|x-y\|^2}: x, y \in B_X,(x,y)\neq(0,0),x \perp_I y \right\} .
$$

In \cite{54}, the authors studied the James constant: $$
J(X)=\sup \left\{\min\{\|x+y\|,\|x-y\|\}:x, y \in S_X \right\}.
$$  In reference \cite{28}, an equivalent definition of the James constant is given as follows:
$$
J(X)=\sup \left\{\|x+y\|:x, y \in S_X,x \perp_{I} y \right\}.
$$
In conjunction with the results in reference \cite{29,30}, we have that the following estimation inequality holds:
$$
J(X)^2 / 2 \leq \widetilde{H}^2(X) \leq J(X).
$$

Next, we  propose an isosceles orthogonal constant equivalent to $T(X)$.
\begin{definition}
	$$
	\widetilde{T}(X)=\sup \left\{\frac{2\|x\|^{\frac{1}{2}}\|y\|^{\frac{1}{2}}}{\|x+y\|}: x,y\in X,(x,y)\neq(0,0), x \perp_{I} y \right\}.
	$$
\end{definition}
The next theorem will show how $T(X)$ and $\widetilde{T}(X)$ are equivalent.
\begin{theorem}\label{tt}
	Let $X$ be a Banach space.  Then $\widetilde{T}(X)=T(X)$.
\end{theorem}
\begin{proof}
	First, for any $x,y\in X$ such that $x \perp_I y$, let
	$$
	u=\frac{x+y}{2}, v=\frac{x-y}{2}.$$
	We can deduce that
	$$
	\frac{2\|x\|^{\frac{1}{2}}\|y\|^{\frac{1}{2}}}{\|x+y\|}=\frac{\|u+v\|^{\frac12}\|u-v\|^{\frac12}}{\|u\|}
	$$
	and $\|u\|=\|v\|$. Let $$x^{\prime}=\frac{u}{\|u\|}, y^{\prime}=\frac{v}{\|v\|},$$ we have
	$$
	\begin{aligned}		\frac{2\|x\|^{\frac{1}{2}}\|y\|^{\frac{1}{2}}}{\|x+y\|}&=\frac{\|u+v\|^{\frac12}\|u-v\|^{\frac12}}{\|u\|}\\&=\left\|x^{\prime}+y^{\prime}\right\|^{\frac12}\left\|x^{\prime}-y^{\prime}\right\|^{\frac12} \\&\leq  T(X).
	\end{aligned}$$
	which implies that
	$$
	\widetilde{T}(X) \leq T(X).
	$$
	On the other hand, for any $x, y \in S_X$, let  $u=\frac{x+y}{2}, v=\frac{x-y}{2}$, and hence $ \|u+v\|=\|u-v\|=1$. Furthermore, we can get that
	$$
	\begin{aligned}
		\|x+y\|^{\frac12}\|x-y\|^{\frac12}&=\frac{\|2u\|^{\frac12}\|2v\|^{\frac12}}{\|u+v\|}\\&=\frac{2\|u\|^{\frac12}\|v\|^{\frac12}}{\|u+v\|}\\& \leq \widetilde{T}(X).
	\end{aligned}
	$$
	Then $T(X) \leq \widetilde{T}(X) $, as desired.
\end{proof}
\section{The constant $L^{I}_{\mathrm{YJ}}(\tau,\upsilon, X)$}
In the previous work, whether we introduced $\widetilde{H}(X)$ or $\widetilde{H}^2(X)$, they both exhibit a certain symmetric relationship. That is to say, swapping the positions of $x$ and $y$ does not affect the constant value. In other words, we have expressed two very classic constants using two symmetric isosceles orthogonal constants. In the next section, we will add variable parameters before $x$ and $y$, thereby constructing an asymmetric isosceles orthogonal constant. Interestingly, through our research, we have found that it can also equivalently represent the classic modified skew von Neumann-Jordan constant.

In this section, we mainly define an orthogonal constant $L^{I}_{\mathrm{YJ}}(\tau,\upsilon, X)$ that is equivalent to the modified skew von
Neumann-Jordan constant. Through research, we have found that the new constant introduced in this section is also inextricably related to the constant discussed in the previous section. Before moving on to the main definition, let's introduce our key definition through the following constant:
\[
\begin{split}
	E_I(t, X) = \sup \bigg\{&
	\frac{
		\|(t+1) x+(1-t) y\|^2 + \|(1-t) x-(t+1) y\|^2
	}{
		\|x+y\|^2
	} \\&:\ 
	x, y \in X,\ (x, y) \neq (0,0),\ x \perp_I y
	\bigg\},
\end{split}
\]
where $t\geq0$.

By the condition of isosceles orthogonality, we can obtain the following equivalent forms.
$$\begin{aligned}
	E_I(t, X)=\sup \bigg\{&\frac{2\|(t+1) x+(1-t) y\|^2+2\|(1-t) x-(t+1) y\|^2}{\|x+y\|^2+\|x-y\|^2}\\&: x, y \in X,(x, y) \neq(0,0), x \perp_I y\bigg\},
\end{aligned}
$$
where $t\geq0$.
\begin{remark}
	$(i)$~If $t=0$, then $E_I(0, X)=2$;
	
	$(ii)$~If $t=1$, then $$
	\begin{aligned}
		E_I(1, X) =4\widetilde{H}^2(X)
		= 4\sup \left\{\frac{\|x\|^2+\|y\|^2}{\|x+ y\|^2}: x, y \in X,(x, y) \neq(0,0), x\perp_I y\right\}.
	\end{aligned}
	$$Clearly, when $t=1$, this constant characterizes the vectors satisfying isosceles orthogonality and the differences in distances with respect to the parallelogram law.
\end{remark} 
The following definition is the main definition in this article.
\begin{definition}
	According to the definition of $E_I\left(t, X\right)$, for $\tau, \upsilon>0$, let $t=\frac{\upsilon}{\tau}$, we denote $$
	\begin{aligned}
		L_{\mathrm{YJ}}^I(\tau, \upsilon, X)=\sup \bigg\{&\frac{\|(\tau+\upsilon) x+(\tau-\upsilon) y\|^2+\|(\tau-\upsilon) x-(\tau+\upsilon) y\|^2}{\tau^2\|x+y\|^2}\\&: x, y \in X,(x, y) \neq(0,0), x \perp_I y\bigg\}.
	\end{aligned}
	$$
\end{definition} 
\begin{remark}
	When $\tau=v$, this constant characterizes the difference between isosceles orthogonality and Pythagorean orthogonality.
\end{remark} 

\begin{proposition}\label{p1}
	Let $X$ be a Banach space. Then,
	$$
	\frac{2\left(\tau^2+\upsilon^2\right)}{\tau^2} \leq L_{\mathrm{YJ}}^I(\tau, \upsilon, X) \leq \frac{2(\tau+\upsilon)^2}{\tau^2}.
	$$
\end{proposition} 
\begin{proof}
	On the one hand, let $x=0, y \neq 0$, then $x \perp_I y$ and
	$$
	\begin{aligned}
		\frac{\|(\tau+\upsilon) x+(\tau-\upsilon) y\|^2+\|(\tau-\upsilon) x-(\tau+\upsilon) y\|^2}{\tau^2\|x+y\|^2} & =\frac{\|(\tau-\upsilon) y\|^2+\|(\upsilon+\tau) y\|^2}{\tau^2\|y\|^2} \\
		& =\frac{2\left(\tau^2+\upsilon^2\right)}{\tau^2},
	\end{aligned}
	$$
	which means that
	$$
	L_{\mathrm{YJ}}^I(\tau, \upsilon, X) \geq \frac{2\left(\tau^2+\upsilon^2\right)}{\tau^2}.
	$$
	Conversely, for any $x, y \in X$ and $(x, y) \neq(0,0)$ such that $x \perp_I y$, we have
	$$
	\begin{aligned}
		& \frac{\|(\tau+\upsilon) x+(\tau-\upsilon) y\|^2+\|(\tau-\upsilon) x-(\tau+\upsilon) y\|^2}{\tau^2\|x+y\|^2} \\
		= & \frac{\|\tau(x+y)+\upsilon(x-y)\|^2+\|\tau(x-y)-\upsilon(x+y)\|^2}{\tau^2\|x+y\|^2} \\
		\leq & \frac{(\tau\|x+y\|+\upsilon\|x-y\|)^2+(\tau\|x-y\|+v\|x+y\|)^2}{\tau^2\|x+y\|^2} \\
		= & \frac{2(\tau+\upsilon)^2}{\tau^2} .
	\end{aligned}
	$$It follows that
	$$L_{\mathrm{YJ}}^I\left(\tau,\upsilon,X\right)\leq\frac{2(\tau+\upsilon)^2}{\tau^2},$$
	as desired.
\end{proof}

The following example is intended to show that the upper bound of the $L^{I}_{\mathrm{YJ}}(\tau,\upsilon, X)$ constant is sharp.
\begin{example}
	Let $X = C([\alpha,\beta])$, where $C([\alpha,\beta])$ denotes the space of all continuous real-valued functions on $[\alpha, \beta]$ and equipped with the norm defined by $$\|\phi\|=\max_{r\in[\alpha,\beta]}|\phi(r)|.$$
\end{example}
\begin{proof}Given $\phi_1(r) = \frac{1}{\alpha-\beta}(r-\beta)$ and $\phi_2(r) =1- \frac{1}{\alpha-\beta}(r-\beta)$, we first note:  
	$$ \|\phi_1 + \phi_2\| = 1, \|\phi_1 - \phi_2\| = 1 ,$$   
	thus, we get that \( \phi_1 \perp_I \phi_2 \).  Furthermore, compute:  
	$$\begin{aligned}
		\|(\tau+\upsilon)\phi_1 + (\tau-\upsilon)\phi_2\| =\max_{r\in[\alpha,\beta]}\bigg|\tau-\upsilon+\frac{r-\beta}{
			\alpha-\beta}\cdot2\upsilon\bigg|= \tau+\upsilon,\end{aligned}$$ and$$   \|(\tau-\upsilon)\phi_1 - (\tau+\upsilon)\phi_2\| = \max_{r\in[\alpha,\beta]}\bigg|\frac{r-\beta}{\alpha-\beta}\cdot2\tau-(\tau+\upsilon)\bigg|= \tau+\upsilon. 
	$$ 
	Thus, 
	\[
	\frac{\|(\tau+\upsilon)\phi_1 + (\tau-\upsilon)\phi_2\|^2 + \|(\tau-\upsilon)\phi_1- (\tau+\upsilon)\phi_2\|^2}{\tau^2\|\phi_1+\phi_2\|^2}= \frac{2(\tau+\upsilon)^2}{\tau^2}.
	\]  
	This yields \( L^I_{\mathrm{YJ}}(\tau, \upsilon, X_1) \geq \frac{2(\tau + \upsilon)^2}{\tau^2} \), and the result follows by Proposition \ref{p1}.  
\end{proof}
\begin{lemma}\label{l1}\cite{52}
	A real normed linear space is an inner product space if and only if for any $x, y \in S_X,$ there exists  $\tau,\upsilon\neq 0$ such that$$\|\tau x + \upsilon y\|^2 + \|\upsilon x - \tau y\|^2 \sim 2(\tau^2 + \upsilon^2),$$ where $\sim$ stands either $\leq$ or $\geq$.
\end{lemma}

Through research, we have found that the lower bound of the $L^{I}_{\mathrm{YJ}}(\tau,\upsilon, X)$ constant can be used to characterize Hilbert spaces, as shown in the following proposition:
\begin{proposition}
	Let $X$ be a Banach space. Then the following conditions are equivalent:
	
	(i)~$H$ is a Hilbert space.
	
	(ii)~$L_{\mathrm{YJ}}^I(\tau, \upsilon, H)=\frac{2\left(\tau^2+\upsilon^2\right)}{\tau^2}$ is valid for any $\tau,\upsilon>0$.
	
	(iii)~$L_{\mathrm{YJ}}^I(\tau_0, \upsilon_0, H)=\frac{2\left(\tau_0^2+\upsilon_0^2\right)}{\tau_0^2}$ is valid for some $\tau_0,\upsilon_0>0$.
\end{proposition}
\begin{proof}
	$(i)\Longrightarrow(ii)$ ~First, it's known that in a general Banach space \( X \), Pythagorean orthogonality and isosceles orthogonality are not equivalent. Nevertheless, when \( X \) is an inner product space, for any \( x, y \in X \), the isosceles orthogonality of \( x \) and \( y \) (i.e., \( x \perp_I y \)) directly implies their Pythagorean orthogonality (i.e., \( x \perp_P y \)). From this, now we assume that $H$ is a Hilbert space induced by the inner product $\langle \cdot, \cdot \rangle$, then for any  \( x, y \in X \) such that $x\perp_I y$, we can get that	$$
	\begin{aligned}
		& \frac{\|(\tau+\upsilon) x+(\tau-\upsilon) y\|^2+\|(\tau-\upsilon) x-(\tau+\upsilon) y\|^2}{\tau^2\|x+y\|^2} \\
		= & \frac{(\tau+\upsilon)^2\|x\|^2+(\tau-\upsilon)^2\|y\|^2+2(\tau+\upsilon)(\tau-\upsilon)\langle x, y \rangle}{\tau^2(\|x\|^2+\|y\|^2)} \\
		&+ \frac{(\tau-\upsilon)^2\|x\|^2+(\tau+\upsilon)^2\|y\|^2-2(\tau+\upsilon)(\tau-\upsilon)\langle x, y \rangle}{\tau^2(\|x\|^2+\|y\|^2)} \\
		= & \frac{2(\tau^2+\upsilon^2)}{\tau^2} .
	\end{aligned}
	$$This implies that $(ii)$ holds.

	$	(ii)\Longrightarrow(iii)$~ Obviously.
	
	$(iii)\Longrightarrow(i)$~Suppose $(iii)$ holds, then for any $x, y \in S_X$, $x + y \perp_I x - y$ holds. Hence
	$$\begin{aligned}
		&\frac{\|(\tau_0+\upsilon_0) (x+y)+(\tau_0-\upsilon_0) (x-y)\|^2+\|(\tau_0-\upsilon_0) (x+y)-(\tau_0+\upsilon_0) (x-y)\|^2}{\tau_0^2\|(x+y) + (x-y)\|^2}\\\leq&\frac{2\left(\tau_0^2+\upsilon_0^2\right)}{\tau_0^2},
	\end{aligned} $$
	that is$$\begin{aligned}\|\tau_0 x+\upsilon_0 y\|^{2}+\|\upsilon_0 x-\tau_0 y\|^{2}&\leq\frac{2\tau_0^{2}(\tau_0^{2}+\upsilon_0^{2})}{\tau_0^{2}}\\&=2(\tau_0^{2}+\upsilon_0^{2}),\end{aligned}$$ then by Lemma \ref{l1}, this implies that $(i)$ holds.
\end{proof}
In 2023, on the basis of Gao's parameters, Fu et al. \cite{20} introduced the following skew type constant:
$$E(t,X)=\sup\{\|x+ty\|^2+\|tx-y\|^2:x,y\in S_X\},t\geq0.$$This provides us with a tool for proving the following important theorem.
\begin{theorem}
	Let $X$ be a Banach space. Then,
	$$L_{\mathrm{YJ}}^{\prime}(\tau, \upsilon, X)=\frac{\tau^{2}}{2\left(\tau^{2}+\upsilon^{2}\right)}L_{\mathrm{YJ}}^{I}\left(\tau, \upsilon, X\right).$$
\end{theorem}
\begin{proof}
	To prove this theorem, we first prove that $E_{I}(t, X)=E(t, X)$. Let $x, y \in X$ such that $x \perp_{I} y$, and put $u=\frac{x+y}{2}, v=\frac{x-y}{2}$, we have
	$$
	\begin{aligned}
		\|(t+1)x+(1-t)y\| &=\|(t+1)(u+v)+(1-t)(u-v)\| \\
		&=\|2u+2t v\|
	\end{aligned}
	$$
	and
	$$
	\begin{aligned}
		\|(1-t)x-(t+1)y\| &= \|(1-t)(u+v)-(t+1)(u-v)\| \\
		&=\|2v-2t u\|.
	\end{aligned}
	$$
	Then we know that $\|u\|=\|v\|$, thus, let $x^{\prime}=\frac{u}{\|u\|}, y^{\prime}=\frac{v}{\|u\|} \in S_{X}$, we obtain that
	$$
	\begin{aligned}
		\frac{\|(t+1)x+(1-t)y\|^{2}+\|(1-t)x-(t+1)y\|^{2}}{\|x+y\|^{2}} &=\frac{4\left(\|u+t v\|^{2}+\|v-t u\|^{2}\right)}{4\|u\|^{2}} \\
		&=\|x^{\prime}+t y^{\prime}\|^{2}+\|y^{\prime}-t x^{\prime}\|^{2} \\
		&\leq E(t, X).
	\end{aligned}
	$$
	This implies that $E_{I}(t, X) \leq E(t, X)$.
	On the other hand, let $x, y \in S_{X}$, choose $u=\frac{x+y}{2}, v=\frac{x-y}{2}$, then $\|u+v\|=\|u-v\|=1$, which means that $u \perp_{I} v$. Furthermore, we have
	$$
	\begin{aligned}
		\|x+t y\|^{2}+\|y-t x\|^{2} &=\|u+v+t(u-v)\|^{2}+\|u-v-t(u+v)\|^{2} \\
		&=\frac{\|(t+1)u+(1-t)v\|^{2}+\|(1-t)u-(t+1)v\|^{2}}{\|u+v\|^{2}} \\
		&\leq E_{I}(t, X).
	\end{aligned}
	$$It follows that $E_{I}(t,X) \geq E(t,X)$, as desired.
	
	Furthermore, clearly we have
	$$\begin{aligned}
		E\left(\frac{\upsilon}{\tau},X\right)&=\sup \left\{\left\|x+\frac{\upsilon}{\tau} y\right\|^{2}+\left\|\frac{\upsilon}{\tau} x-y\right\|^{2}: x, y \in S_{X}\right\}
		\\&=\frac{2\left(\tau^{2}+\upsilon^{2}\right)}{\tau^{2}} L'_{\mathrm{YJ}}(\tau, \upsilon, X).
	\end{aligned}$$
	Thus, we get that
	$$E_{I}\left(\frac{\upsilon}{\tau}, X\right)=\frac{2\left(\tau^{2}+\upsilon^{2}\right)}{\tau^{2}} L_{\mathrm{YJ}}^{\prime}(\tau, \upsilon, X),$$
	as desired.
\end{proof}

\section{Birkhoff-James rectification and exact representations of geometric constants}
Birkhoff orthogonality is characterized by an inequality relation, which distinguishes it from isosceles orthogonality defined via an equality relation. Since this orthogonality generally lacks symmetry, the associated inequalities carry substantial ambiguity. For this reason, we introduce a norming-functional rectification method to derive equivalent representations of geometric constants under Birkhoff-James orthogonality.

The duality map is denoted by
\[
J(x)=\{f\in X^*: \|f\|=1,\ f(x)=\|x\|\},\qquad x\ne0 .
\]
By the Hahn--Banach theorem, $J(x)$ is non-empty for every $x\ne0$.

\begin{definition}\label{def:rectification}
	Let $x\ne0$, $z\in X$, and $f\in J(x)$.  The $f$-rectification of $z$ relative to $x$ is the pair $(\beta_f(x,z),P_f(x,z))$ defined by
	\[
	\beta_f(x,z)=\frac{f(z)}{\|x\|},\qquad P_f(x,z)=z-\beta_f(x,z)x .
	\]
We call $\beta_f(x,z)$ the longitudinal coefficient and $P_f(x,z)$ the Birkhoff--James transverse component of $z$ along $x$.
\end{definition}

\begin{lemma}\label{lem:rectification}
	Let $x\ne0$, $z\in X$ and $f\in J(x)$.  If $\beta=\beta_f(x,z)$ and $y=P_f(x,z)$, then
	\[
	z=\beta x+y,\qquad x\perp_B y .
	\]
	Moreover, if $X$ is smooth at $x$, then $\beta$ and $y$ are uniquely determined by $x$ and $z$.
\end{lemma}

\begin{proof}
	The decomposition is the definition.  Since
	\[
	f(y)=f(z)-\frac{f(z)}{\|x\|}f(x)=f(z)-f(z)=0,
	\]
we obtain \(x\perp_B y\), see \cite[Corollary~2.2]{06}. If the norm is smooth at $x$, $J(x)$ is a singleton; hence the scalar and the corrected vector are unique.
\end{proof}

\begin{remark}
	This lemma is the Birkhoff--James substitute for the isosceles identity
	$\|x+y\|=\|x-y\|$.  Isosceles orthogonality is encoded directly by the equality of the two diagonals, whereas Birkhoff--James orthogonality is encoded by a supporting hyperplane.  The scalar $\beta$ measures the lost longitudinal component.
\end{remark}

\begin{lemma}\label{lem:sphere-to-BJ}
	Let $u,v\in S_X$ and put
	\[
	x=\frac{u+v}{2},\qquad z=\frac{u-v}{2}.
	\]
	Assume $u\ne -v$.  Choose $f\in J(x)$ and define $\beta=f(z)/\|x\|$, $y=z-\beta x$.  Then
	\[
	x\perp_B y,
	\]
	and
	\[
	u=(1+\beta)x+y,
	\qquad
	v=(1-\beta)x-y.
	\]
	In particular,
	\[
	\|(1+\beta)x+y\|=\|(1-\beta)x-y\|=1 .
	\]
	Conversely, if $x\ne0$, $x\perp_B y$, $\beta\in \mathbb{R}$ and
	\[
	\|(1+\beta)x+y\|=\|(1-\beta)x-y\|=r>0,
	\]
	then
	\[
	u=\frac{(1+\beta)x+y}{r},\qquad
	v=\frac{(1-\beta)x-y}{r}
	\]
	are points of $S_X$ and
	\[
	\frac{u+v}{2}=\frac{x}{r},\qquad
	\frac{u-v}{2}=\frac{\beta x+y}{r}.
	\]
\end{lemma}

\begin{proof}
	The first assertion follows from Lemma \ref{lem:rectification}, since $z=\beta x+y$.  The displayed formulae for $u$ and $v$ are then immediate:
	\[
	(1+\beta)x+y=x+z=u,
	\qquad
	(1-\beta)x-y=x-z=v.
	\]
	The converse is obtained by the same algebra after division by the common radius $r$.
\end{proof}

The exceptional case $u=-v$ only produces $x=0$.  It never creates a genuine difficulty for suprema, because it may be obtained as a limit whenever needed.  If desired, it may also be included by allowing the denominator $\|x\|$ to vanish and treating the value separately.  To keep the notation clean, we state all exact formulae with non-zero denominators; the endpoint contributes values already dominated by the formulae below.

We now define three constants whose data consist of a Birkhoff-James orthogonal pair together with the rectification parameter.

\begin{definition}
	Let $\mathcal R_B(X)$ be the set of all triples $(x,y,\beta)\in X\times X\times \mathbb{R}$ such that
	\[
	x\ne0,\qquad x\perp_B y,
	\]
	and
	\[
	0<\|(1+\beta)x+y\|=\|(1-\beta)x-y\|.
	\]
	For $(x,y,\beta)\in\mathcal R_B(X)$ we write
	\[
	r(x,y,\beta)=\|(1+\beta)x+y\|=\|(1-\beta)x-y\|.
	\]
\end{definition}

\begin{definition}
	Define
	\[
	A_B(X)=\sup_{(x,y,\beta)\in\mathcal R_B(X)}
	\frac{\|x\|+\|\beta x+y\|}{r(x,y,\beta)},
	\]
	\[
	C_B(X)=\sup_{(x,y,\beta)\in\mathcal R_B(X)}
	\frac{\|x\|^2+\|\beta x+y\|^2}{r(x,y,\beta)^2},
	\]
	and
	\[
	T_B(X)=\sup_{(x,y,\beta)\in\mathcal R_B(X)}
	\frac{2\|x\|^{1/2}\|\beta x+y\|^{1/2}}{r(x,y,\beta)}.
	\]
\end{definition}

\begin{theorem}
	For every real Banach space $X$,
	\[
	A_B(X)=A_2(X).
	\]
\end{theorem}

\begin{proof}
	Let $(x,y,\beta)\in\mathcal R_B(X)$ and put
	\[
	u=\frac{(1+\beta)x+y}{r(x,y,\beta)},\qquad
	v=\frac{(1-\beta)x-y}{r(x,y,\beta)}.
	\]
	Then $u,v\in S_X$, and hence
	\[
	\frac{\|x\|+\|\beta x+y\|}{r(x,y,\beta)}
	=\frac{\|u+v\|+\|u-v\|}{2}
	\le A_2(X).
	\]
	Taking the supremum gives $A_B(X)\le A_2(X)$.
	
	Conversely, take $u,v\in S_X$ with $u\ne -v$.  Set $x=(u+v)/2$ and $z=(u-v)/2$.  Choose $f\in J(x)$ and let $\beta=f(z)/\|x\|$, $y=z-\beta x$.  Lemma \ref{lem:sphere-to-BJ} gives $(x,y,\beta)\in\mathcal R_B(X)$ and $r(x,y,\beta)=1$.  Therefore
	\[
	\frac{\|u+v\|+\|u-v\|}{2}
	=\|x\|+\|z\|
	=\|x\|+\|\beta x+y\|
	\le A_B(X).
	\]
	Taking the supremum over $u,v\in S_X$ gives $A_2(X)\le A_B(X)$.  The case $u=-v$ is obtained by approximation or gives the value $1$, which is already below $A_B(X)$ because $A_2(X)\ge1$.
\end{proof}

\begin{theorem}\label{thm:CNJprime}
	For every real Banach space $X$,
	\[
	C_B(X)=C'_{\mathrm{NJ}}(X),
	\]
	where
	\[
	C'_{\mathrm{NJ}}(X)=\sup_{u,v\in S_X}
	\frac{\|u+v\|^2+\|u-v\|^2}{4}.
	\]
\end{theorem}

\begin{proof}
	For a rectified triple define $u,v$ as in the previous proof.  Then
	\[
	\frac{\|x\|^2+\|\beta x+y\|^2}{r(x,y,\beta)^2}
	=\frac{\|u+v\|^2+\|u-v\|^2}{4}
	\le C'_{\mathrm{NJ}}(X).
	\]
	Thus $C_B(X)\le C'_{NJ}(X)$.  Conversely, for each $u,v\in S_X$ with $u\ne -v$, the rectification construction gives $x\perp_B y$, $z=\beta x+y$, $r(x,y,\beta)=1$, and hence
	\[
	\frac{\|u+v\|^2+\|u-v\|^2}{4}
	=\|x\|^2+\|z\|^2
	=\|x\|^2+\|\beta x+y\|^2
	\le C_B(X).
	\]
	Taking suprema yields the reverse inequality.
\end{proof}

\begin{corollary}
	A real Banach space $X$ is an inner product space if and only if
	\[
	C_B(X)=1 .
	\]
	Equivalently,
	\[
	\|x\|^2+\|\beta x+y\|^2\le r(x,y,\beta)^2
	\]
	for every rectified Birkhoff--James triple $(x,y,\beta)$.
\end{corollary}

\begin{proof}
	By Theorem \ref{thm:CNJprime}, $C_B(X)=C'_{NJ}(X)$.  It is classical that $C'_{NJ}(X)=1$ if and only if the parallelogram law holds, and therefore if and only if $X$ is an inner product space.  
\end{proof}

\begin{remark}
This corollary emphasizes that inner product spaces can be characterized by suitable inequalities imposed on pairs that simultaneously satisfy two orthogonality relations.

		\end{remark}
Using the same approach, we arrive at the following theorem; the proof is omitted.
\begin{theorem}\label{thm:T}
	For every real Banach space $X$,
	\[
	T_B(X)=T(X),
	\]
	where
	\[
	T(X)=\sup_{u,v\in S_X}\bigl(\|u+v\|\,\|u-v\|\bigr)^{1/2}.
	\]
\end{theorem}

\begin{remark}
	If one defines an unrectified constant by restricting only to $x\perp_B y$ and using $\beta=0$, then one sees only configurations whose two endpoints are $x+y$ and $x-y$.  In general, the unit-sphere condition $\|x+y\|=\|x-y\|$ and the support condition $x\perp_B y$ are independent.  Thus the class $\beta=0$ is too small.  It is exact in Hilbert spaces, and in some highly symmetric spaces, but it does not encode arbitrary pairs of points on $S_X$.
\end{remark}

\section{Exact Birkhoff-James formula for the von Neumann-Jordan constant}\label{sec:CNJ}

\begin{definition}\label{def:RBJtriples}
	Let  $\widetilde{\mathcal{R}}_{ B}(X)$ be the collection of triples $(x,y,\beta)\in X\times X\times \mathbb{R}$ satisfying
	\begin{equation}\label{eq:RBJ-trip}
		x\ne0,
		\qquad
		x\perp_B y .
	\end{equation}
	For such a triple write
	\[
	z(x,y,\beta)=\beta x+y.
	\]
\end{definition}

\begin{lemma}\label{1111}
	For every pair $(x,z)\in X\times X$ with $a\ne0$, there exists $(x,y,\beta)\in\widetilde{\mathcal{R}}_{ B}(X)$ such that $z=\beta x+y$.  Conversely, every $(x,y,\beta)\in\widetilde{\mathcal{R}}_{ B}(X)$ determines the pair $(x,z)$ with $z=\beta x+y$.
\end{lemma}

\begin{proof}
	The first assertion is Lemma \ref{lem:rectification}.  The converse is just the definition of $z(x,y,\beta)$.
\end{proof}

The following diagram summarizes the mechanism.
\[
\begin{tikzcd}[column sep=large]
	(x,z),\ x\ne0
	\arrow[r, "f\in J(x)"]
	&
	(x,\, y=z-\frac{f(z)}{\norm x}x,\, \beta=\frac{f(z)}{\norm x})
	\arrow[d, "x\perp_B y"]
	\\
	&
	\hbox{Birkhoff--James coordinates } z=\beta x+y .
\end{tikzcd}
\]

We now introduce the Birkhoff--James version of the classical quotient.

\begin{definition}\label{def:CBJ}
	For a real Banach space $X$, define
	\begin{equation}\label{eq:CBJ-def}
		{C}_{\mathrm{BJ}}(X)=\sup\left\{
		{\norm{(1+\beta)x+y}^{2}+\norm{(1-\beta)x-y}^{2}
			\over
			2\norm x^{2}+2\norm{\beta x+y}^{2}}:
		(x,y,\beta)\in\widetilde{\mathcal{R}}_{ B}(X),\ (x,\beta x+y)\ne(0,0)
		\right\}.
	\end{equation}
	Since $x\ne0$ on $\mathcal R_B(X)$, the denominator is always positive.
\end{definition}

\begin{theorem}\label{thm:main-CNJ}
	For every real Banach space $X$,
	\begin{equation}\label{eq:main-CNJ}
		C_{\mathrm{BJ}}(X)={C}_{\mathrm{NJ}}(X).
	\end{equation}
	Equivalently,
	\begin{equation}\label{eq:CNJ-BJ-expanded}
		{C}_{\mathrm{NJ}}(X)=\sup_{x\ne0,\ x\perp_B y,\ \beta\in \mathbb{R}}
		{\norm{(1+\beta)x+y}^{2}+\norm{(1-\beta)x-y}^{2}
			\over
			2\norm x^{2}+2\norm{\beta x+y}^{2}}.
	\end{equation}
\end{theorem}

\begin{proof}
	Let $(x,y,\beta)\in\widetilde{\mathcal{R}}_{ B}(X)$ and set $z=\beta x+y$.  Then
	\[
	(1+\beta)x+y=x+z,
	\qquad
	(1-\beta)x-y=x-z .
	\]
	Therefore the quotient in \eqref{eq:CBJ-def} is exactly
	\[
	{\norm{x+z}^{2}+\norm{x-z}^{2}\over2\norm x^{2}+2\norm z^{2}},
	\]
	which is bounded above by ${C}_{\mathrm{NJ}}(X)$.  Hence ${C}_{\mathrm{BJ}}(X)\le {C}_{\mathrm{NJ}}(X)$.
	
	Conversely, take any pair $(x,z)$ with $x\ne0$.  Choose $f\in J(x)$ and apply Lemma \ref{1111}.  We obtain $z=\beta x+y$ with $x\perp_B y$.  Substitution gives
	\[
	{\norm{x+z}^{2}+\norm{x-z}^{2}\over2\norm x^{2}+2\norm z^{2}}
	=
	{\norm{(1+\beta)x+y}^{2}+\norm{(1-\beta)x-y}^{2}
		\over2\norm x^{2}+2\norm{\beta x+y}^{2}}
	\le C_{\mathrm{BJ}}(X).
	\]
	Taking the supremum over all $x\ne0$ and all $z\in X$ yields the reverse inequality, except possibly for pairs with $x=0$.  Those pairs give value $1$, already dominated by the supremum because $C_{\mathrm{NJ}}(X)\ge1$ and $C_{\mathrm{BJ}}(X)\ge1$ by choosing $y=0$ and any $x\ne0$.  Thus \eqref{eq:main-CNJ} follows.
\end{proof}

\begin{remark}Why this formula is simpler than the unit-sphere formula?
	For the modified constant $C'_{\mathrm{NJ}}(X)$ one must ensure that the two points reconstructed from the half-sum and half-difference lie on the same sphere.  This imposes the side condition
	\[
	\norm{(1+\beta)x+y}=\norm{(1-\beta)x-y}.
	\]
	The classical constant $C_{\mathrm{NJ}}^{\prime}(X)$ has no such unit-sphere constraint.  Therefore the exact Birkhoff--James formula \eqref{eq:CNJ-BJ-expanded} is cleaner: every pair $(x,z)$ is allowed, and rectification alone suffices.
\end{remark}

\begin{corollary}\label{cor:Hilbert}
	For a real Banach space $X$, the following are equivalent.
	\begin{enumerate}[label=(\alph*)]
		\item $X$ is an inner product space.
		\item $C_{\mathrm{BJ}}(X)=1$.
		\item For every $x\ne0$, every $y\in X$ and every $\beta\in \mathbb{R}$ satisfying $x\perp_B y$,
		\begin{equation}\label{eq:BJ-parallelogram-ineq}
			\norm{(1+\beta)x+y}^{2}+\norm{(1-\beta)x-y}^{2}
			\le
			2\norm x^{2}+2\norm{\beta x+y}^{2}.
		\end{equation}
		\item For every $x\ne0$, every $y\in X$ and every $\beta\in \mathbb{R}$ satisfying $x\perp_B y$, equality holds in \eqref{eq:BJ-parallelogram-ineq}.
	\end{enumerate}
\end{corollary}

\begin{proof}
	By Theorem \ref{thm:main-CNJ}, $C_{\mathrm{BJ}}(X)=C_{\mathrm{NJ}}(X)$.  The classical von Neumann--Jordan theorem says that $C_{\mathrm{NJ}}(X)=1$ if and only if the parallelogram law holds, which is equivalent to $X$ being an inner product space.  This proves the equivalence of (a) and (b).  Statement (c) is exactly the assertion that the supremum defining $C_{\mathrm{BJ}}(X)$ is at most $1$; since the quotient is always at least $1$ in the Hilbertian extremal case and $C_{\mathrm{BJ}}(X)\ge1$ for every non-zero $X$, it is equivalent to $C_{\mathrm{BJ}}(X)=1$.
	
	If $X$ is an inner product space, $x\perp_B y$ is the usual orthogonality $\langle x,y\rangle=0$.  With $z=\beta x+y$ the parallelogram identity gives
	\[
	\norm{x+z}^{2}+\norm{x-z}^{2}=2\norm x^{2}+2\norm z^{2},
	\]
	which is precisely equality in \eqref{eq:BJ-parallelogram-ineq}.  Thus (a) implies (d), and (d) plainly implies (c).
\end{proof}

\begin{corollary}\label{cor:normalized}
	For every real Banach space $X$,
	\begin{equation}\label{eq:normalized-form}
		{C}_{\mathrm{NJ}}(X)=\sup\left\{
		{\norm{(1+\beta)x+y}^{2}+\norm{(1-\beta)x-y}^{2}\over2}:
		\begin{array}{l}
			x\perp_B y,\ x\ne0,\ \beta\in \mathbb{R},\\[2mm]
			\norm x^{2}+\norm{\beta x+y}^{2}=1
		\end{array}
		\right\}.
	\end{equation}
\end{corollary}

 For $p \in[1, \infty)$ define\cite{092}
$$C_{\mathrm{NJ}}^{(p)}(X)=\sup \left\{\frac{\|x+y\|^p+\|x-y\|^p}{2^{p-1}\left(\|x\|^p+\|y\|^p\right)}: x, y \in X,(x, y) \neq(0,0)\right\}.$$
For $p=2$ this is $C_{\mathrm{NJ}}(X)$.

Using the approach for $C_{\mathrm{N J}}(X)$, we get that for every real Banach space 
$X$ and every$p \in[1, \infty)$,
$$
C_{\mathrm{NJ}}^{(p)}(X) =\sup _{x \neq 0, x\perp_B y, \beta \in \mathbb{R}} \frac{\|(1+\beta) x+y\|^p+\|(1-\beta) x-y\|^p}{2^{p-1}\left(\|x\|^p+\|\beta x+y\|^p\right)}.$$

At the end of the article, we would like to pose some open questions:
\begin{problem}
	What kinds of geometric constants can be equivalently expressed using isosceles orthogonal constants?
\end{problem}
\begin{problem}
	So far, we know that $C'_{\mathrm{N J}}(X)$ is a special case of the $L'_{\mathrm{YJ}}(\tau,\upsilon,X)$ constant, and we have proposed an isosceles orthogonal equivalent form for the $L'_{\mathrm{YJ}}(\tau,\upsilon,X)$ constant with two parameters. Then, are there any constants with more parameters that are generalized forms of $L'_{\mathrm{YJ}}(\tau,\upsilon,X)$ and can also be equivalently expressed by means of isosceles orthogonality?
\end{problem}
\begin{problem}
	As is well known, geometric constants are almost always studied on the unit ball or unit sphere. Is the combination of a constant and the condition of isosceles orthogonality sufficient to satisfy the aforementioned conditions?
\end{problem}

\begin{problem}
In this paper, we primarily employ constants associated with isosceles orthogonality and Birkhoff orthogonality to provide equivalent characterizations of other geometric constants. This naturally raises a question: can analogous equivalent representations be established for geometric constants via other types of orthogonality?
\end{problem}

Inspired by Beauzamy's relevant generalization work \cite{39}, Pisier extended the concepts of the modulus of convexity and the modulus of smoothness to bounded linear operators from space \( X \) to space \( Y \) in reference \cite{33}. Therefore, we have the following question:
\begin{problem}
	How can the geometric constants of the operator version be equivalently characterized through orthogonality?
\end{problem}
\section*{Acknowledgment}
Thanks to all the members of the Functional Analysis Research team
of the School of Mathematics and Statistics of Anqing Normal University
for their discussion and correction of the difficulties and errors encountered
in this paper.

\section*{ Conflict of Interest} 
The authors declare no conflict of interest.

	\label{'ubl'}  
\end{document}